%
%
%
%
   \input amstex
\documentstyle{amsppt}
\magnification=\magstep1
\parindent=1em
\CenteredTagsOnSplits
\NoBlackBoxes
\nopagenumbers
\NoRunningHeads
\pageno=1
\footline={\hss\tenrm\folio\hss}
\def\lp{\bold l^{\bold p}}
        \topmatter
\title {On factorization of operators through the spaces $ \lp$ }
\endtitle
       \author { Oleg I. Reinov${{ }^\dag}$}  \endauthor
\address\newline
Oleg I. Reinov \newline
Department of Mathematics\newline
St Petersburg University\newline
St Peterhof, Bibliotech pl 2\newline
198904  St Petersburg, Russia
\endaddress

\email
orein\@orein.usr.pu.ru
\endemail

\thanks
${{ }^\dag}$This work was done with partial support of the Ministry of the
general and professional education of Russia (Grant 97-0-1.7-36) and
FCP ``Integracija", reg. No. 326.53.
\endthanks

\abstract
We give conditions on a pair of Banach spaces
$ X$ and $ Y,$  under which
each operator from $ X$ to $ Y,$ whose second adjoint
factors compactly through the space $ l^p,\,$
$ 1\le p\le+\infty,\,$ itself compactly factors through $ l^p.$
The conditions are as follows:
either the space $ X^*,$ or the space $ Y^{***}$ possesses
the Grothendieck approximation property.
Leaving the corresponding question for parameters $ p>1,\,p\neq 2,$
still open, we show that for $ p=1$
the conditions are essential.
\endabstract
    \endtopmatter

\document
\baselineskip=18pt
\footnote""{${ }^\ddag$
AMS Subject Classification:  47B10. Hilbert--Schmidt operators,
trace class operators, nuclear operators, p-summing operators, etc.}

\footnote""{${ }$
Key words: $p$-compact operators,
bases, approximation property, tensor products.}

\magnification=\magstep1
\voffset1ex

\def\a{\alpha}
\def\ffi{\varphi}
\define\e{\varepsilon}
 \def\ot{\otimes}
 \def\wh{\widehat}
 \def\wt{\widetilde}
 \def\tto{\wt{\wt\ot}}
\def\({\left(} \def\){\right)} \def\[{\left[} \def\]{\right]}
   \def\<{\langle}
   \def\>{\rangle}

\def\tr{\operatorname{trace}\,}
\def\N{\frak N}
\def\K{\frak K}
\def\L{\frak L}
\def\D{\frak D}
\def\bw{\text{\bf w}}
\def\bK{\text{\bf K}}
\def\bN{\text{\bf N}}

\bigskip

One of few questions, which remain open in the theory
of operator ideals and are connected with the notion of their regularity,
is the question of the possibility of a factorization
of a continuous operator, acting between Banach spaces,
through the classical sequence spaces $ l^p$ under the condition
that this operator, considered as an operator with values
in the second dual to the space of images, admits a factorization
of such a kind.
More precisely, the question is connected with
the possibility of the "compact factorization" through these spaces
(and we will deal with this "compact" situation).
Clearly, we reject the trivial case where $ p=2,$
since in this case the answer is affirmative by evident reasons:
every (closed) subspace of a Hilbert space is complemented.

A few of known facts concerning the formulated question
is two results obtained by the author yet in 1982 (see [1],
Corollaries 3.4 and 4.2):
if $ p=1$  or $ p=+\infty$ (the factorization through the spaces
$ l^1$ and $ c_0)$ then the answer is negative.
The technique developed in [1] allowed us to construct the corresponding
counterexamples in the spaces which, in general,
do not possess the Grothendieck approximation property.

In this note we shall give sufficient (and likely close to necessary)
conditions, imposed on Banach spaces under consideration,
under which the answer to the mentioned question is positive
for each value of the parameter $p,\, 1\le p\le \infty.$
On the other hand we will show that the conditions which appear bellow
can not be rejected (or weakened). But, unfortunately,
the counterexample that we will give concerns only one case, ---
the case of the compact factorization through the space $ l^1.$

As before, {\it the question on the $l^p$-factorization of operators
in the case where $ p>1,\, p\neq 2,$ still}
(for a long time) {\it remains open}.

It seems to me that I have a counterexample,
similar to that of Theorem 2 bellow,
for the case of the factorization through the space
$ c_0,$ but since the corresponding proof is not completely checked up
yet, we do not consider it in the note.
\smallpagebreak

We hold standard notation of the geometrical theory of operators
in Banach spaces. The classical reference book on the theory of operator
ideals is the A. Pietsch's monograph [2].
All the spaces under consideration are Banach spaces.
We usually denote the elements of the spaces by the corresponding small
letters: $x\in X,\,y\in Y,\dots,\, x'\in X^*,\,y''\in Y^{**},\dots$.
For  $ p\in [1,+\infty],$
the conjugate exponent  $ p'$ is defined by the relation $ 1/p+1/p'=1.$
By $\L\left(X,Y \right)$ it is denoted the space of all
(linear continuous) operators from $X$ to $Y$ with its standard norm
(generally, in this paper the term "operator" always means
"linear continuous operator").

Recall that {\it a Banach space has the approximation property}
\, (of A. Grothendieck)
if the identity map on the space can be approximated,
in the topology of compact convergence, by finite dimensional
operators (some details and other reformulations of this definition
can be found either in the book [2] or directly in the origin [3]).
And one more important reminder. If $ J$ is an operator ideal,
then $ J^{ \operatorname{ reg}}(X, Y)$ denotes
the space of all operators $ T$ from $ X$ to $ Y,$ for which
$ \pi_Y\,T\in J(X, Y^{**}),$ where $ \pi_Y$ is
the canonical isometric imbedding of the space $ Y$
into its second dual $ Y^{**}.$
In some cases it will be convenient to consider some of the Banach
spaces as the subspaces of their second duals
(under the canonical imbeddings),
without mentioning the operators $ \pi$ which give these imbeddings.
\smallpagebreak

Let $ 0<p\le\infty.$ A family $ \left\{ x_j\right\}_{j=1}^\infty,$
where $ x_j\in X,$  is said to be {\it weakly $ p$-summable}
if $ \left\{ \<x_j, x' \>\right\}\in l^p$ for all $ x'\in X^*.$
Following A. Pietsch and his notations, we set
$$ \bw_p(x_j):=
  \sup_{\|x'\|\le 1} \left(\sum_j \, |\<x_j,\,x'\>|^p\right)^{1/p}.
$$

If $ 1\le p\le\infty$ and
$ S$ is a finite dimensional operator from $ X$ to $ Y,$
then its $ \K_p^0$-norm is the quantity
$\bK_p^0(S):= \inf \left\{ \sup_i |\a_i|\,\bw_p(x'_j)\,\bw_{p'}(y_j)
\right\},$
where the infimum is taken over all {\it finite} representations
of the operator $ S$ in the form
$$ S=\sum_{j=1}^n \a_j\,x'_j\ot y_j.
$$
We shall denote by $ X^*\tto_p Y$
the completion of the algebraic tensor product $ X^*\ot Y$
(which we consider as a linear space of all finite dimensional
operators) in this norm $ \bK_p^0,$
keeping the same notation $ \bK_p^0$
for the norm in this completion.
Each tensor element $ z\in  X^*\tto_p Y$ can be represented in the form
of a converging series
$ z=\sum_{j=1}^\infty \a_j\,x'_j\ot y_j,$ where $ \a_j\to 0,$
and, moreover, for a given
$ \e>0$ this representation can be chosen in such a way that
$$  \sup_i |a_i|\,\bw_p(x'_j)\,\bw_{p'}(y_j) \le \bK_p^0(z)+\e.
$$

The tensor product $ X^*\tto_p Y$ generates naturally
a linear subspace of operators in $ \L(X,Y);$
the operator corresponding to an element
$ z$ will be denoted by $ \tilde z.$
It is quite possible (but, in general, it is not known whether it is),
that some of
{\it different tensors} \, $ z_1, z_2$ may be "pasting together"
when turning into the operators
(i. e. for $ z_1\neq z_2$ it is not forbidden for the equality
$ \tilde z_1 = \tilde z_2$ to be valid).
We denote (again, following the notation of A. Pietsch) by
$ \K_p(X,Y)\,$ the factor space
of the considered tensor product over the kernel
(if non trivial) of the mentioned mapping and
by $ \bK_p$ the obtained norm on this space of operators.
Note that $ \[\K_p, \bK_p\] $ is a normed operator ideal
and it is a particular case of the so called ideals
of $ (r, p ,q)$-nuclear operators, ---
$ \[\N_{(r,p,q)}, \bN_{(r,p,q)}\]\!:$\ \,
$ \[\K_p, \bK_p\] =  \[\N_{(\infty,p,p')}, \bN_{(\infty,p,p')}\].$
For details we refer the reader to the book [2].
However, there is a thing that is very important for us and
which is concerned with the possibility of the mentioned
nontrivial factorization. Namely,
{\it if one of the space $ X^*$ or $ Y$ has the approximation
property then the following {\rm (formal)} equality holds}:\,
$ \[X^*\tto_p Y, \bK_p^0\] = \[\K_p(X,Y), \bK_p\],$
that is the kernel of the canonical map
$ X^*\tto_p Y \to \K_p(X, Y)$ is trivial and this map is
an isometric bijection. The proof of this fact is now standard,
so we omit it suggesting to the reader as a not difficult exercise.

We say that an operator $ T: X\to Y$ {\it compactly factors
through the space $ l^p$} \, ({\it $p$-compact}\ in the terminology of
A. Pietsch), if there exist two compact operators
$ A: X\to l^p$ and $ B: l^p\to Y$ for which $ T=BA.$ As a norm
of such an "$l^p$-factored operator" we take the quantity
$ \inf\,\|A\|\,\|B\|,$  where the infimum is taken over all possible
factorizations of the operator $ T$ of the mentioned kind.

One of the central theorem in the theory of compactly
$ l^p$-factored operators is formulated with the help of the following
relations:
$$ \bK_p(T) = \inf\,\{\|A\|\,\|B\|:\, T=BA: X\to l^p\to Y;\, A,B \,
\text{ are compact} \}
$$
and {\it the ideal $ \K_p$ coincides with the ideal of all
operators which can be compactly factored through
}\ $ l^p$\, (see [2], Theorem 18.3.2).
The space $ \D_{p'}(Y,X^{**})\, $ which is conjugate to
$ X^*\tto_p Y$ is described in the same place of the mentioned
monograph, but since, as a matter of fact,
we need not to know of what operators this
conjugate space consists, we are not going to give more details.
It is worthwhile, however, for the reader
to look in the book for getting a minimal information
about the elements of the conjugate space.
For us it is important to know that from the definition
of the operators lying in this space it follows immediately that
{\it all the series that will appear in the proof of our first theorem
are convergent}\, (so we won't have to be distracted
every time to check up the convergence of each appearing series).

\proclaim {Theorem {\rm 1}}
Let $\, p\in [1,+\infty].$
If one of the space $\, X^* $ or $\, Y^{***} $ has the
Grothendieck approximation property, then any operator
from $ X$  into $ Y,$ whose second adjoint factors compactly
through the space $ l^p,$ itself compactly factors through
\,$ l^p.$ Thus, in this case
$ \K^{ \operatorname{ reg}}_p (X, Y)= \K_p (X, Y).$
\endproclaim

\demo{Proof}
   Suppose that there exists such an operator
$ T\in \L(X,Y)$ that
$ T\notin \K_p(X,Y),$ but $ \pi_Y\,T\in \K_p(X,Y^{**}).$
Since either $ X^*$ or $ Y^{**}$ has the approximation property,
then $ \K_p(X,Y^{**})=X^*\tto_p Y^{**}.$
Therefore, the operator $ \pi_Y\,T$ can be identified with
a tensor element $ t\in X^*\tto_p Y^{**};$ in addition,
by the choice of $ T,$ \
$  t\notin X^*\tto_p Y$ \ (here $  X^*\tto_p Y$ is considered as
a subspace of the space $  X^*\tto_p Y^{**}$).
Hence, there exists such an operator
$ U\in \D_{p'}(Y^{**},X^{**})=\( X^*\tto_p Y^{**}\)^*$ that
$ \tr U\circ t=\tr \(t^*\circ \( U^*|_{X^*}\) \)=1$ and
$ \tr U\circ \pi_Y\circ z=0$ for every $ z\in X^*\tto_p Y.$
>From the last it follows, in particular, that
$$ U|_{\pi_Y(Y)}=0 \, \text{ É }\,  \pi_Y^*\,U^*|_{X^*}=0. \tag1$$
Indeed, if $ x'\in X^*$ and $ y\in Y,$ then
$$ \<U\pi_Y\,y,x'\> = \<y, \pi_Y^*\,U^*|_{X^*}x'\> =
\tr \,(U\pi_Y)\circ (x'\otimes y)=0.
$$
Evidently, the tensor element $ U\circ t$  generates the
operator $ U\pi_Y T$ which is, by the previous reasoning, equal
identically to zero.

If the space $ X^*$ has the approximation property then
 $ X^*\tto_p X^{**}= \K_p(X,X^{**})$
and so the tensor element $ U\circ t$ is zero, in contradiction
with the equality $ \tr\, U\circ t=1.$

Let now $ Y^{***}$ has the approximation property.
In this case
$ V:= \( U^*|_{X^*}\)\circ T^*\circ \pi_Y^*:$ \
$$ \CD    Y^{***} @>\pi_Y^* >> Y^* @> T^* >> X^* @> U^*|_{X^*} >> Y^{***}
\endCD $$
uniquely determines a tensor element
$ t_0\in Y^{****}\tto_{p'} Y^{***}.$
Let us take any representation $ t=\sum \a_n x'_n\otimes y''_n$ for $ t$
as an element of the space $ X^*\tto_p Y^{**}.$
Denoting for the simplicity by $ U_*$ the operator $ U^*|_{X^*}$ and
recalling that $ \pi_Y T=t,$ we get:
$$\multline
   Vy'''=U_*\, \( T^*\pi_Y^*\,y'''\) =
    U_*\, \( (T^*\pi_Y^*\pi_{Y^*})\,\pi_Y^*\, y'''\) = \\ =
    U_*\, \( (\pi_Y T)^*\,\pi_{Y^*})\,\pi_Y^*\, y'''\) =
    U_*\, \( \(\sum \a_n y''_n\otimes x'_n\)\,\pi_{Y^*}\,\pi_Y^*\, y'''\)
  = \\ =U_*\, \( \sum \a_n \<y''_n, \pi_Y^*\, y'''\> \,x'_n \) =
      \sum \a_n \<\pi_Y^{**}y''_n,  y'''\> \,U_* x'_n.
  \endmultline
$$

Thus, the operator $ V$ (or the element $ t_0$) has in the space
$ Y^{****}\tto_{p'} Y^{***}$ the representation
$$ V= \sum \a_n \pi_Y^{**}(y''_n)\otimes U_* (x'_n).
$$
Therefore (see (1)),
$$  \tr t_0=\tr V= \sum \a_n \<\pi_Y^{**}(y''_n), U_* (x'_n)\> =
        \sum \a_n \<y''_n, \pi_Y^*\,U_* x'_n\> =  \sum 0=0.
$$
On the other hand,
$$\multline
 Vy'''= U_* \( \pi_Y T\)^* y'''= U_*\circ t^* (y''')=\\
     =U_*\, \( \sum \a_n \<y''_n, y'''\> \, x'_n\)=
      \sum \a_n \<y''_n, y'''\> \, U_* x'_n,
\endmultline $$
whence $ V=\sum y''_n\otimes U_*(x'_n).$  Therefore
$$ \tr t_0=\tr V= \sum \a_n \<y''_n, U_* x'_n\> =
 \sum \a_n \<Uy''_n, x'_n\> = \tr U\circ t=1.$$
The obtained contradiction completes the proof of the theorem.
 $\quad\blacksquare$ \enddemo

Our next result shows that, in general, the approximation
conditions imposed on $ X$ and $Y\,$
in Theorem 1 are essential and can not be replaced
by weaker conditions "one of the space
$ X$ or $ Y^{**}$ has the (bounded) approximation property".

\proclaim {Theorem {\rm 2}}
There exists a Banach space $Z,\,$ with the following properties:
    \newline    \phantom{iii} i\,{\rm)}
all the duals spaces of $ Z$ are separable;
    \newline    \phantom{ii} ii\,{\rm)}
$Z^{**}\,$ has a boundedly complete basis;
    \newline    \phantom{i} iii\,{\rm)}
$Z^{***}\,$ does not have the approximation property;
    \newline    \phantom{i} iv\,{\rm)}
there exists an operator $ T\in\L(Z^{**},Z)$
that is not factored compactly through the space $ l^1,$
and for which, however, the operator
$ \pi_Z\,T: Z^{**}\to Z\hookrightarrow Z^{**}$ compactly factors
through $ l^1.$ \newline\phantom{iii}
 Thus,
    $ \K_1(Z^{**}, Z)\subsetneqq  \K_1^{ \operatorname{ reg}}(Z^{**}, Z).$
\endproclaim
\demo{Proof}
We use the following result obtained in the paper [1]
(see Lemma 1.2,2) and Corollary 1.2) in [1]):
there exists a reflexive separable Banach space
$ E$ such that for every $ r,\,$ $ 1\le r\le \infty\,$ $ r\neq 2,$\,
the canonical map $ E^*\wh\ot_r E\to \N_r(E, E)$ is not one-to-one.
Here $ E^*\wh\ot_r E$ is the $ r$-projective tensor product
associated with the space $ \N_r(E, E)$ of all
$ r$-nuclear operators in $ E.$ We consider the case $ r=\infty$ and,
passing to the conjugate spaces, formulate the dual result
(setting $ X:=E^*$):
there exists a reflexive separable Banach space
$ X$ such that the canonical mapping $ X^*\tto_1 X\to \K_1(X, X)$
is not one-to-one.

By using the theorem 1 and the corollary 1 from the paper [4]
(see also the proof of the mentioned corollary),
find such a separable Banach space $ Z$
that $ Z^{**}$ has a boundedly complete basis, the space $ X$
is isomorphic to the factor space $ Z^{**}/\pi_Z(Z)$
and, in addition, the corresponding linear homomorphism
$ \ffi: Z^{**}\to X$ \, (with the kernel
$ \operatorname{ Ker \ffi}= \pi_Z(Z)$) has the property that
the subspace $ \ffi^*(X^*)$ is complemented in $ Z^{***}.$

Let $ u\in  X^*\tto_1 X$ be a nonzero tensor element generating
an operator $ \wt u$ which is identically zero.
Let us consider any representation of the tensor
$ u$ in the form
$$ u= \sum_{j=1}^\infty \a_j\,x'_j\ot x_j,
$$
where the sequence $ ( x_j)$ is bounded,
the sequence $ ( x'_j)$
is weakly 1-summable and $ \a_j\to 0.$
We can lift up the elements $ x_j$ to the elements
$ z''_j$ of the space $ Z^{**}$ in such a way that for every
$ j\in\Bbb N$\ $ \ffi(z''_j)=x_j$ and the sequence
$ (z''_j)$ is bounded in $ Z^{**}.$
Arising from $ u$ a new tensor element
$$ v = \sum_{j=1}^\infty \a_j\,x'_j\ot z''_j,
$$
belongs to the tensor product $ X^*\tto_1 Z^{**}$
and is not equal to zero,
and also the composition $ u=\ffi\circ v,$ considered
as an operator, is zero.
Since the space $ Z^{**}$ has a basis (and, in particular,
has the approximation property), the element $ v$ can be correctly
identified with the associated operator $ \wt v: X\to Z^{**},$
which is thus not zero and for which, as said above, the image
$ \wt v(X)$ lies in the subspace $ \pi_Z(Z),$ the kernel of the
operator $ \ffi.$
Thus, we get an operator $ V: X\to Z$ for which
$V=\pi^{-1}_Z \wt v_0 $ and $ V^{**}=\wt v,$
where $ \wt v_0: X\to \pi_Z(Z)$ is the operator induced by
$ \wt v.$
It is not difficult to see that the operator $ V$
can not belong to the space $ \K_1(X, Z)=X^*\tto_1 Z,$
since else, for each its (nonzero!) tensor representation
in this space, its image under acting of the homomorphism
$ \ffi$ would be a zero tensor in $ X^*\tto_1 X$
(we identify now the spaces $ Z$ and $ \pi_Z(Z)$);
and on the other hand
this image is just $ u.$
So, we can conclude that
$X^*\tto_1 Z = \K_1(X, Z)\subsetneqq \K_1^{ \operatorname{ reg}}(X, Z). $

Finally, recalling that the subspace $ \ffi^*(X^*)$ is complemented
in $ Z^{***},$ we see that $ \K_1(Z^{**}, Z) =
Z^{***}\tto_1 Z\subsetneqq \K_1^{ \operatorname{ reg}}(Z^{**}, Z).$
 $\quad\blacksquare$ \enddemo
\newpage

\centerline{REFERENCES}
\smallpagebreak

\ref \no 1
\by O. I. Reinov \pages   125-134
\paper   Approximation properties of order p and the existence of
  non-p-nuclear operators with p-nuclear second adjoints
\yr 1982 \vol  109
\jour    Math. Nachr.
\endref

\ref\no 2
\by A. Pietsch
\inbook  Operator ideals
\publaddr North-Holland
\publ
\yr  1980
\endref

\ref \no3
\by A. Grothendieck \pages 196+140
\paper  Produits tensoriels topologiques et espases nucl\'eaires
\yr 1955\vol  16
\jour  Mem. Amer. Math. Soc.
\endref

\ref \no 4
\by L. Lindenstrauss \pages  279-284
\paper  On James' paper "Separable Conjugate Spaces"
\yr 1971\vol 9
\jour Israel J. Math.
\endref

\enddocument